\documentclass[a4paper,12pt]{article}
\usepackage{amssymb,amsmath,amsthm, amscd}
\usepackage{verbatim}

\DeclareMathOperator{\Spec}{Spec} \DeclareMathOperator{\Aut}{Aut} \DeclareMathOperator{\PGL}{PGL}
\DeclareMathOperator{\GL}{GL} \DeclareMathOperator{\PSL}{PSL} \DeclareMathOperator{\Char}{char}
\newcommand{\Z}{{\mathbb Z}}
\newcommand{\R}{{\mathbb R}}
\newcommand{\C}{{\mathbb C}}
\newcommand{\F}{{\mathbb F}}
\newcommand{\Kbar}{{\overline{K}}}

\DeclareMathOperator{\Gal}{Gal}

\newcommand{\el}[1]{{}^{#1\negmedspace}}

    \newtheorem{thm}{Theorem}
    \numberwithin{thm}{section}
    \newtheorem{cor}[thm]{Corollary} 
    \newtheorem{lem}[thm]{Lemma}
    \newtheorem{prop}[thm]{Proposition}

    \theoremstyle{definition}
    
    \theoremstyle{remark}

    \newcommand{\ol}{\overline}

 \usepackage{sectsty}
\allsectionsfont{\centering\large}
\newenvironment{mat2}{\left(\begin{array}{cc}}{\end{array}\right)}

\begin{document}

\title{\large\textbf{{Fields of Moduli of Hyperelliptic Curves}}}
\author{\small{Bonnie Huggins}}
\date{}

 \maketitle
\begin{abstract}Let $F$ be an algebraically closed field with $\Char(F)\ne 2$, let $F/K$ be a Galois extension, and
  let $X$ be a hyperelliptic curve defined over $F$.  Let $\iota$
be the hyperelliptic involution of
  $X$.  We show that $X$ can be defined over its field of moduli relative to the extension $F/K$
  if $\Aut(X)/\langle \iota\rangle$ is not cyclic.  We
construct explicit examples of hyperelliptic curves not
  definable over their field of moduli when $\Aut(X)/\langle \iota\rangle$ is cyclic.
\end{abstract}

\section{Introduction}
\indent\par Let $X$ be a curve of genus $g$ defined over a field $F$, let $F/L$ be a Galois extension, and let
$K$ be the field of moduli relative the the extension $F/L$.  (See Section 2 for the definition of ``field of
moduli''.) It is well known that if $g$ is $0$ or $1$ then $X$ admits a model defined over $K$.  It is also well
known that if the group of automorphism of $X$ is trivial then $X$ can be defined over $K$; for example, see
Example~1.7 in~\cite{kock:1}. However, if $g\ge2$ and $|\Aut(X)|>1$, the curve $X$ may not be definable over its
field of moduli.

We examine the case where $X$ is hyperelliptic and $F$ is an algebraically closed field of characteristic not
equal to $2$. In this case $\Aut(X)$ is always nontrivial since it contains the hyperelliptic involution
$\iota$. Examples of hyperelliptic curves not definable over their field of moduli are given on page~177
in~\cite{sh:1}. In~\cite{Cardona:1} it is shown that $X$ can be defined over $K$ if $g=2$ and $|\Aut(X)|>2$. In
Theorem~4.2 and Corollary~4.4 of~\cite{Shaska:1} it is shown that $X$ is definable over $K$ if $\Char(F)=0$,
$g\ge 2$, and $\Aut(X)/\langle \iota\rangle$ has at least two involutions.  In Section~1 of~\cite{Shaska:1} it
is conjectured that $X$ is definable over $K$ if $\Char(F)=0$ and $|\Aut(X)|>2$. In this paper, we refute this
conjecture and show that $X$ can be defined over $K$ if $\Aut(X)/\langle \iota\rangle$ is not a cyclic group.
\section{Fields of Moduli}
\indent\par Let $K$ be a field, let $F/K$ be a Galois extension and let $X$ be a hyperelliptic curve defined
over $F$. Let $\sigma\in\Gal(F/K)$.  The curve $\el{\sigma}{X}$ is the base extension $X \underset{\Spec
F}\times \Spec F$ of
 $X$ by the morphism $\Spec F\xrightarrow{\Spec\sigma}\Spec F$.
The field of moduli relative to the extension $F/K$ is defined as the fixed field $F^{H}$ of
\[H:=\{\sigma\in \Gal(F/K)\mid X\cong \el{\sigma}{X}\mbox{ over } F\}.\]
A subfield $E$ of $F$ is a field of definition for $X$ if there exists a curve $X_E$ defined over $E$ such that
$X\cong X_E \underset{\Spec E}\times \Spec F$.

\begin{prop}\label{closed}  Let $K_m$ be the field of moduli of $X$.  Then the
subgroup $H$ is a closed subgroup of $\Gal(F/K)$ for the Krull topology.  That is,
\[H=\Gal(F/K_m).\]
 The field of $K_m$ is contained in
each field of definition between $K$ and $F$ (in particular, $K_m$ is a finite extension of $K$).  Hence if the
field of moduli is a field of definition, it is the smallest field of definition between $F$ and $K$.  Finally,
the field of moduli of $X$ relative to the extension $F/K_m$ is $K_m$.

\end{prop}
\begin{proof}
See Proposition 2.1 in~\cite{Debes:1}.
\end{proof}

\section{Finite Subgroups of 2-Dimensional Projective General Linear Groups}
\indent\par Throughout this section let $\ol K$ be an algebraically closed field of characteristic $p$ with
$p=0$ or
 $p>2$.  In the following two lemmas we identify a matrix in $\GL_2(\ol K)$ with its image in $\PGL_2(\ol K)$.
\begin{lem}\label{subgroups}
Any finite subgroup $G$ of $\PGL_2(\overline{K})$ is conjugate to one of the following groups:

\begin{enumerate}
\item[Case I:] when $p=0$ or $|G|$ is relatively prime to $p$.
\begin{enumerate}

\item[(a)]  $G_{C_n}=\Big{\{}\begin{mat2}

  \zeta^r & 0 \\
0 & 1 \\
\end{mat2}: r=0,1,\ldots,n-1\Big{\}}\cong C_n,~n\ge 1$

\item [(b)]$G_{D_{2n}}=\Big{\{}\begin{mat2}

  \zeta^r & 0 \\
0 & 1 \\
\end{mat2},
\begin{mat2}

  0 & \zeta^r \\
1 & 0 \\
\end{mat2}: r=0,1,\ldots,n-1\Big{\}}\cong D_{2n},~{n>2}$

\item [(c)] $G_{V_4}=\Big{\{}\begin{mat2}
\pm 1 & 0 \\
0 & 1 \\
\end{mat2},
\begin{mat2}
0 & \pm 1 \\
1 & 0 \\
\end{mat2}\Big{\}}\cong V_4:= \Z/2\Z \times \Z/2\Z$

\item [(d)] $G_{A_4}=\Big{\{}\begin{mat2}
\pm 1 & 0 \\
0 & 1 \\
\end{mat2},
\begin{mat2}
0 & \pm 1 \\
1 & 0 \\
\end{mat2},
\begin{mat2}
  i^{\nu} & i^{\nu} \\
1 & -1 \\
\end{mat2},
\begin{mat2}
i^{\nu} & -i^{\nu} \\
1 & 1 \\
\end{mat2},$

$\begin{mat2}
1 & i^{\nu} \\
1 & -i^{\nu} \\
\end{mat2},
\begin{mat2}
-1 & -i^{\nu} \\
1 & -i^{\nu} \\
\end{mat2}:\nu=1,3\Big{\}}\cong A_4$

\item [(e)] $G_{S_4}=\Big{\{}\begin{mat2}
  i^{\nu} & 0 \\
0 & 1 \\
\end{mat2},
\begin{mat2}

  0 & i^{\nu} \\
1 & 0 \\
\end{mat2},
\begin{mat2}
i^{\nu} & -i^{\nu+\nu'} \\
1 &  i^{\nu'} \\
\end{mat2}: \nu, \nu'=0,1,2,3\Big{\}}{\;\cong S_4}$

\item [(f)] $G_{A_5}=\Big{\{}\begin{mat2}
  \epsilon^r & 0 \\
0 & 1 \\
\end{mat2},
\begin{mat2}

  0 & \epsilon^r \\
-1 & 0 \\
\end{mat2},
\begin{mat2}
  \epsilon^r\omega & \epsilon^{r-s} \\
1 & -\epsilon^{-s}\omega \\
\end{mat2},
\begin{mat2}
  \epsilon^r\overline{\omega} & \epsilon^{r-s} \\
1 & -\epsilon^{-s}\overline{\omega} \\
\end{mat2}: r, s = 0,1,2,3,4\Big{\}}\cong A_5$
\end{enumerate}
\noindent where $\omega=\frac{-1+\sqrt{5}}{2}$, $\overline{\omega}=\frac{-1-\sqrt{5}}{2}$, $\zeta$ is a
primitive $n^{th}$ root of unity, $\epsilon$ is a primitive $5^{th}$ root of unity, and $i$ is a primitive
$4^{th}$ root of unity.
 \item [Case II:] when $|G|$ is divisible by $p$.
\begin{enumerate}
  \item[(g)]
$G_{\beta,A}=\left\{ \begin{mat2}
                      \beta^k &a\\
                      0&1\end{mat2}\colon a\in A, \; k\in\Z\right\}$, where $A$ is a finite
                      additive subgroup of $\ol K$ containing $1$ and $\beta$ is a root of unity such that $\beta A=A$
\item[(h)] $\PSL_2(\mathbb{F}_{p^r})$
  \item[(i)] $\PGL_2(\mathbb{F}_{p^r})$
\end{enumerate}
where $\mathbb{F}_{p^r}$ is the finite field with $p^r$ elements.
\end{enumerate}
\end{lem}
\begin{proof}
See \textsection\textsection 55-58 in~\cite{Weber:1} and Chapter 3 in~\cite{Suzuki:1}.
\end{proof}
\begin{lem}\label{normalizer} Let $N(G)$ be the normalizer of $G$ in
$\PGL_2(\overline{K})$.  Then
\begin{enumerate}
\item [(a)] $N(G_{C_n})=\Big\{\begin{mat2}
\alpha & 0 \\
0 & 1 \\
\end{mat2},
\begin{mat2}
0 & \alpha \\
1 & 0 \\
\end{mat2}: \alpha\in \overline{K}^*\Big\}$ if $n>1$,
\item[(b)] $N(G_{D_{2n}})=G_{D_{4n}}$ if $n>2$, \item[(c)] $N(G_{V_4})=G_{S_4}$, \item[(d)]$N(G_{A_4})=G_{S_4}$,
\item[(e)] $N(G_{S_4})=G_{S_4}$, \item[(f)] $N(G_{A_5})=G_{A_5}$,
\item[(h)]$N(\PSL_2(\mathbb{F}_{p^r}))=\PGL_2(\mathbb{F}_{p^r}),\mbox{ and }$ \item[(i)]
$N(\PGL_2(\mathbb{F}_{p^r}))=\PGL_2(\mathbb{F}_{p^r})$.
\end{enumerate}
\end{lem}
\begin{proof} \hfill
\begin{enumerate}
\item[(a)] See \textsection 55 in~\cite{Weber:1}. \item[(b)] See \textsection 55 in~\cite{Weber:1}. \item[(c)]
 Since $G_{V_4}$ is a normal subgroup of $G_{S_4}$, $G_{S_4}\subseteq N(G_{V_4})$. Conjugation of $G_{V_4}$ by
$G_{S_4}$ gives a homomorphism $G_{S_4}\to\Aut(V_4)\cong S_3$.  A computation shows that the centralizer $Z$ of
$G_{V_4}$ in $\PGL_2(\Kbar)$ is $G_{V_4}$.  The kernel of this homomorphism is $Z\cap G_{S_4}=Z$.  Since
$G_{S_4}/Z\cong S_3$, every automorphism of $G_{V_4}$ is given by conjugation by an element of $G_{S_4}$. Let
$U\in N(G_{V_4})$. Then $UV\in Z=G_{V_4}$ for some $V\in G_{S_4}$, so $U\in G_{S_4}$. \item[(d)] Since $G_{V_4}$
is a characteristic subgroup of $G_{A_4}$, $N(G_{A_4})\subseteq N(G_{V_4}){=G_{S_4}}$. As $G_{A_4}$ is normal in
$G_{S_4}$, we get $N(G_{A_4})=G_{S_4}$. \item[(e)] Since $G_{A_4}$ is a characteristic subgroup of $G_{S_4}$,
$N(G_{S_4})\subseteq N(G_{A_4}){=G_{S_4}}$. Thus $N(G_{S_4})=G_{S_4}$. \item[(f)] Conjugation of $G_{A_5}$ by
$N(G_{A_5})$ gives a homomorphism $N(G_{A_5})\to\Aut(A_5)$. The kernel of this homomorphism is the centralizer
of $G_{A_5}$ in $N(G_{A_5})$, which is just the centralizer $Z$ of $G_{A_5}$ in $\PGL_2(\Kbar)$. A computation
shows that $Z$ is just the identity. Since $\Aut(A_5)$ is finite, $N(G_{A_5})$ is a finite subgroup of
$\PGL_2(\overline{K})$.  Since $G_{A_5}\subseteq N(G_{A_5})$, by Lemma~\ref{subgroups} we must have
$N(G_{A_5})=G_{A_5}$.

\item[(h)] We first show that $N(\PSL_2(\mathbb{F}_{p^r}))$ is finite.  Conjugation of
$\PSL_2(\mathbb{F}_{p^r})$ by $N(\PSL_2(\mathbb{F}_{p^r}))$ gives a homomorphism
$N(\PSL_2(\mathbb{F}_{p^r}))\to\Aut(\PSL_2(\mathbb{F}_{p^r}))$. The kernel of this homomorphism is the
centralizer $Z$ of $\PSL_2(\mathbb{F}_{p^r})$ in $\PGL_2(\Kbar)$.  A computation shows that $Z$ is just the
identity.  Since $\Aut(\PSL_2(\mathbb{F}_{p^r}))$ is finite, so is $N(\PSL_2(\mathbb{F}_{p^r}))$.  By
Lemma~\ref{subgroups} any finite subgroup of $\PGL_2(\Kbar)$ containing $\PSL_2(\mathbb{F}_{p^r})$ must be
isomorphic to either $\PGL_2(\mathbb{F}_{q}))$ or $\PSL_2(\mathbb{F}_{q})$ for some $q$.  Since
$SL_2(\mathbb{F}_{p^r})$ is normal in $\GL_2(\mathbb{F}_{p^r})$, $\PSL_2(\mathbb{F}_{p^r})$ is a normal subgroup
of $\PGL_2(\mathbb{F}_{p^r})$.  So $\PGL_2(\mathbb{F}_{p^r})\subseteq N(\PSL_2(\mathbb{F}_{p^r}))$, in
particular $\PSL_2(\mathbb{F}_{p^r})$ is strictly contained in $N(\PSL_2(\mathbb{F}_{p^r}))$.  By the corollary
on page 80 of~\cite{Suzuki:1}, $\PSL_2(\mathbb{F}_{q})$ is simple for $q>3$.  It follows that
$N(\PSL_2(\mathbb{F}_{p^r}))\ne \PSL_2(\mathbb{F}_{q})$ for $q\ge 3$.  By Theorem 9.9 on page 78
of~\cite{Suzuki:1}, the only nontrivial normal subgroup of $\PGL_2(\mathbb{F}_{q})$ is $\PSL_2(\mathbb{F}_{q})$
if $q>3$.  Therefore $N(\PSL_2(\mathbb{F}_{p^r}))=\PGL_2(\mathbb{F}_{p^r})$. \item[(i)] Clear from the proof of
the previous case.
\end{enumerate}
\end{proof}

\section{Isomorphisms of Hyperelliptic Curves}

\indent\par Throughout this section let $K$ be a perfect field of characteristic $p$ with $p=0$ or $p>2$ and let
$X$ be a hyperelliptic curve defined over an algebraic closure $\overline{K}$ of $K$ with $K$ as its field of
moduli.  In particular, $X$ admits a degree-$2$ morphism to ${\mathbb P}^1$ and the genus of $X$ is at least
$2$.  Each element of Aut($X$)
  induces an automorphism of $\mathbb{P}^1$ fixing the branch points.
The number of branch points is $\ge 3$ (in fact $\ge 6$), so $\Aut(X)$ is finite.  We get a homomorphism
$\Aut(X) \to \Aut({\mathbb P}^1) = \PGL_2(\Kbar)$ with kernel generated by the hyperelliptic involution $\iota$.
Let $G\subset PGL_2({\overline{K}})$ be the image of this homomorphism.  Replacing the original map $X \to
{\mathbb P}^1$ by its composition with an automorphism $g \in \Aut({\mathbb P}^1) = \PGL_2(\Kbar)$ has the
effect of changing $G$ to $gGg^{-1}$, so we may assume that $G$ is one of the groups listed in
Lemma~\ref{subgroups}.
  Fix an equation $y^2=f(x)$ for $X$
  where $f\in \overline{K}[x]$ and $\mbox{disc}(f)\ne 0$.  So the function field $\Kbar(X)$ equals $\Kbar(x,y)$.
\begin{prop}\label{isom}  Let $X$ be as above and let $X'$ be a
hyperelliptic curve defined over $\overline{K}$ given by $y^2=f'(x)$, where $f'(x)$ is another squarefree
polynomial in $\Kbar[x]$.  Every isomorphism $\varphi\colon X\to X'$ is given by an expression of the form:

\[ (x,y)\mapsto \left(\frac{ax+b}{cx+d},
\frac{ey}{(cx+d)^{g+1}}\right),\]

\noindent for some $M=\begin{mat2}

  a & b \\
c & d \\
\end{mat2}\in \GL_2(\overline K)$ and $e\in\overline{K}^*$.
The pair $(M,e)$ is unique up to replacement by $(\lambda M, e \lambda^{g+1})$ for $\lambda\in\overline K^*$. If
$\varphi'\colon X'\to X''$ is another isomorphism, given by $(M',e')$, then the composition $\varphi'\varphi$ is
given by $(M'M,e'e)$.

\end{prop}
\begin{proof}
See Proposition 2.1 in~\cite{Poonen:1}.
\end{proof}

Let $\Gamma=\Gal(\overline{K}/K)$ and let $\sigma\in\Gamma$. Then $\el{\sigma}{X}$ is  the smooth projective
model of $y^2=f^{\sigma}(x)$, where $f^{\sigma}(x)$ is the polynomial obtained from $f(x)$ by applying $\sigma$
to the coefficients.

\begin{lem}\label{isomnormalizer}  Following the notation used above, let $\sigma\in\Gamma$ and suppose that $\varphi\colon X\to
\el{\sigma}{X}$ is given by $(M,e)$.  Let $\overline{M}$ be the image of $M$ in $\PGL_2(\overline{K})$.  If
$G\ne G_{\beta, A}$ then $\overline{M}$ is in the normalizer $N(G)$ of $G$ in $\PGL_2(\overline{K})$.  If
$G=G_{\beta, A}$ then $M$ is an upper triangular matrix.
\end{lem}

\begin{proof}
Since Aut$(\el{\sigma}{X})=\{\psi^{\sigma}\mid\psi\in\Aut(X)\}$, the group of automorphisms of $\mathbb{P}^1$
induced by Aut($\el{\sigma}{X})$ is $G^{\sigma}:=\{U^{\sigma}\mid U\in G\}$.

Let $\psi$ be an automorphism of $X$ given by $(V,v)$. Since $\psi$ is an automorphism, $V\in
\GL_2(\overline{K})$ is a lift of some element $\ol V\in G$. Then $\varphi\psi\varphi^{-1}$ is an automorphism
of $\el{\sigma}{X}$ given by $(MVM^{-1}, v)$.  We have $\ol{MVM^{-1}}{=\ol M\;\ol V\;\ol{M}^{-1}}{\in
G^{\sigma}}$. It follows that $\overline{M}G\overline{M}^{-1}=G^{\sigma}.$   If $G\ne G_{\beta, A}$, by
Lemma~\ref{subgroups}, $G^\sigma=G$. So $\overline{M} \in N(G)$.  If $G=G_{\beta, A}$, then since $G^{\sigma}$
has an elementary abelian subgroup of the same form as $G$, a simple computation shows that $M$ is an upper
triangular matrix.
\end{proof}

\begin{lem}\label{ingsig}
Following the above notation, suppose that for every $\tau\in\Gamma$ there exists an isomorphism
$\varphi_{\tau}\colon X\to \el{\tau}{X}$ given by $(M_{\tau},e)$ where $\overline{M}_{\tau}\in G^{\tau}$. Then
$X$ can be defined over $K$.
\end{lem}
\begin{proof}  Let $P_1,\ldots,P_n$ be the hyperelliptic branch points of $X\to\mathbb{ P}^1$.  Let ${\tau\in\Gamma}$.
The isomorphism $\varphi_{\tau}\colon X\to \el{\tau}{X}$ induces an isomorphism on the canonical images
$\mathbb{P}^1\to\mathbb{P}^1$ which is given by $\overline{M}_{\tau}$.  Then $\ol M_{\tau}$ sends $\{P_1,\ldots,
P_n\}$ to $\{\tau(P_1), \ldots, \tau(P_n)\}$.  Since $\overline{M}_{\tau}\in G^{\tau}$ it merely permutes the
set $\{P_1,\dots,P_n\}$.  Since $\tau$ is arbitrary we have
\[\prod_{P_i\ne\infty } (x-P_i)\in K[x].\]
 It follows that $X$ can be defined over $K$.
\end{proof}
\begin{cor}\label{def}Suppose that $N(G)=G$ and $G\ne G_{\beta, A}$. Then $X$ can be defined over $K$.
\end{cor}
\begin{proof}By Lemma~\ref{subgroups}, $G^{\sigma}=G$ for all $\sigma\in\Gamma$.  Let $\tau\in\Gamma$.  By Lemma~\ref{isomnormalizer}, any isomorphism $X\to\el{\tau}{X}$ is
given by $(M,e)$ where $\ol M\in N(G)=G=G^{\tau}$.

\end{proof}
\section{The Main Result}
\indent\par The following two results of D\`ebes and Emsalem will be used in the proof of our main result.
 They rely on the notions of a cover and the field of moduli of a cover,
for which we refer the reader to \textsection~2.4 in~\cite{Debes:2}.
\begin{thm}\label{bk} Let $F/K$ be a Galois extension and $X$ be a hyperelliptic curve defined over $F$ with $K$ as field
of moduli.  Let $B=X/\Aut(X)$.  Then there exists a model $B_K$ of the curve $B=X/\Aut(X)$ defined over $K$ such
that the cover $X\to B$ with $K$-base $B_K$ is of field of moduli $K$.
\end{thm}
\begin{proof}See Theorem 3.1 in~\cite{Debes:1}.  The authors make the additional assumption that
$\Char(K)$ does not divide $|\Aut(X)|$ but do not use it in their proof.
\end{proof}
\begin{cor}Suppose that $F$ is algebraically closed.  If $B_K$ has a $K$-rational point, then $K$ is a field of definition of $X$.
\end{cor}
\begin{proof} It suffices to show that the cover $X\to B$ with $K$-base $B_K$ can be
defined over $K$, since a field of definition of the cover is automatically a field of definition of $X$. By
Theorem~\ref{bk}, the field of moduli of the cover $X\to B$ with $K$-base $B_K$ is K.  If $K$ is a finite field
then $\Gal(F/K)$ is a projective profinite group.   In this case, by Corollary 3.3 of~\cite{Debes:2} the cover
$X\to B$ can be defined over $K$. If $K$ is not a finite field then since $B_K\cong_K \mathbb {P}^1_K$, $B_K$
has a rational point off the branch point set of $X\to B_K\times F$.  Then by Corollary 3.4 and \textsection~2.9
of~\cite{Debes:2}, the cover can be defined over $K$.
\end{proof}
The curve $B_K$ is called the canonical model of $X/\Aut(X)$ over the field of moduli of $X$. Let
$\Gamma=$Gal($F/K$).  In the proof of Theorem~\ref{bk}, D\`{e}bes and Emsalem show  the canonical model exists
by using the following argument.  For all $\sigma\in\Gamma$ there exists an isomorphism $\varphi_{\sigma}\colon
X\rightarrow \el{\sigma}{X}$ defined over $F$.  Each induces an isomorphism $\tilde{\varphi_{\sigma}}\colon
X/\Aut(X)\rightarrow \el{\sigma}{X}/\Aut(\el{\sigma}{X})$ that makes the following diagram commute:
\[
\begin{CD}
{X} @> {\varphi_{\sigma}} >> {\el{\sigma}{X}} \\
@ V{\rho} VV @ VV {\rho^{\sigma}} V \\
{X/\Aut(X)} @>> {\tilde{\varphi_{\sigma}}} >
{\el{\sigma}{X}/\Aut(\el{\sigma}{X})} \\
\end{CD}
\]
Composing $\tilde{\varphi_{\sigma}}$ with the canonical isomorphism
\[i_{\sigma}\colon \el{\sigma}{X}/\Aut(\el{\sigma}{X})\rightarrow
\el{\sigma}{(X/\Aut(X))}\] we obtain an isomorphism
\[\overline{\varphi_{\sigma}}\colon  X/\Aut(X)\rightarrow
\el{\sigma}{(X/\Aut(X))}.\] The family $\{\overline{\varphi_{\tau}}\}_{\tau\in\Gamma}$ satisfy Weil's cocycle
condition $\overline{\varphi_{\tau}}^{\sigma}\,\overline{\varphi_{\sigma}}=\overline{\varphi_{\sigma\tau}}$
given in Theorem~1 of~\cite{Weil:1}. This shows that $B_K$ exists.

Let $F(B)$ be the function field of $B$. Since $B\cong\mathbb{P}^1$, $F(B)=F(t)$ for some element $t$.  We use
$t$ as a coordinate on $B$.  Suppose that $\overline{\varphi_{\sigma}}$ is given by
\[t\mapsto \frac{at+b}{ct+d}.\] Define $\sigma^*\in\Aut(F(t)/K )$ by
\[\sigma^*(t)=\frac{at+b}{ct+d},\; {\sigma^*(\alpha)}=\sigma(\alpha),\; \alpha\in F.\]  One can verify that
$(\sigma\tau)^*(w)=\sigma^*(\tau^*(w))$ for all $w\in F(t)$.  So we get a homomorphism $\Gamma\to\Aut(F(B)/K)$,
$\sigma\mapsto\sigma^*$. The curve $B_K$ is the variety over $K$ corresponding to the fixed field of
$\Gamma^*=\{\sigma^*\}_{\sigma\in\Gamma}$.

The following lemma will be used in the proof of the main theorem.
\begin{lem}\label{odd}Let $L/K$ be a field extension of odd degree.  Let $C$ be
a curve of genus $0$ defined over $K$ and suppose that $C(L)\ne\varnothing$. Then $C(K)\ne\varnothing$.
\end{lem}

\begin{proof}Let $P\in C(L)$ and let $n=[L:K]$.  Let $\tau_1, \ldots, \tau_n$ be the distinct embeddings of $L$ into an algebraic closure of $L$.  Then
$D=\Sigma\tau_i(P)$ is a divisor of degree $n$ defined over $K$.  Let $\omega$ be a canonical divisor on $C$.
Since deg$(\omega)=-2$, we can take a linear combination of $D$ and $\omega$ to obtain a divisor $D'$ of degree
$1$. Since deg$(\omega-D')<0$, by the Riemann-Roch theorem $l(D')>0$.  So there exists an effective divisor
$D''$ linearly equivalent to $D'$ defined over $K$.  Since $D''$ is effective and of degree 1 it consists of a
point in $C(K)$.
\end{proof}

\begin{thm}\label{main} Let $K$ be a field of characteristic not equal to $2$.
Let $X$ be a hyperelliptic curve defined over $\overline{K}$, an algebraic closure of $K$.  Let
$G=$Aut$(X)/\langle \iota\rangle$ where $\iota$ is the hyperelliptic involution of $X$.  Suppose that $G$ is not
cyclic or that $G$ is cyclic of order divisible by the characteristic of $K$.  Then $X$ can be defined over its
field of moduli relative to the extension $\overline{K}/K$.

\end{thm}
\begin{proof}
Let $\Gamma=\Gal(\overline{K}/K)$.  By Proposition~\ref{closed} we may assume that $K$ is the field of moduli of
$X$. By Proposition~\ref{isom} we may assume that $G$ is given by one of the groups in Lemma~\ref{subgroups}.
Fix an equation $y^2=f(x)$ for $X$
  where $f\in \overline{K}[x]$ and $\mbox{disc}(f)\ne 0$.  So the function field $\Kbar(X)$ equals $\Kbar(x,y)$.
There are eight cases.
\begin{enumerate}
\item[(b)] $G\cong D_{2n}$, $n>2$.  The function field of $X/\Aut(X)$ equals the subfield of $\Kbar(X)$ fixed by
$G_{D_{2n}}$ acting by fractional linear transformations. Then $t:=x^n+x^{-n}$ is fixed by $G_{D_{2n}}$ and is a
rational function of degree $2n$ in $x$, so the function field of $X/\Aut(X)$ equals $\Kbar(t)$. Therefore we
use $t$ as coordinate on $X/\Aut(X)$.  The map $\rho\colon X\rightarrow X/\Aut(X)$ is given by $(x,y)\mapsto
(x^n+x^{-n})$. Let $\sigma\in\Gamma$.  By Lemmas~\ref{isomnormalizer} and~\ref{normalizer},
$\varphi_{\sigma}\colon X\to \el{\sigma}{X}$ is given by $(M,e)$ where $\overline{M}\in D_{4n}$.  Then the map
$\rho^{\sigma}\varphi_{\sigma}\colon X\rightarrow \el{\sigma}{X}/\Aut(\el{\sigma}{X})$ is given by $(x,y)\mapsto
\pm (x^n+x^{-n})$.  So $\sigma^*(t)=\pm t$.  The curve $B_K$ corresponds to the fixed field of $\overline{K}(t)$
under $\Gamma^*$. Then $t=0$ corresponds to a point $P\in B_K(K)$. \item[(c)]  $G\cong V_4$.  The element
$t:=x^2+x^{-2}$ is fixed by $G_{V_4}$ and is a rational function of degree $4$ in $x$.  So the function field of
$X/\Aut(X)$ equals $\Kbar(t)$.  We use $t$ as a coordinate on
 $X/\Aut(X)$.  The map $\rho\colon X\rightarrow X/\Aut(X)$ is given by $(x,y)\mapsto (x^2+x^{-2})$.  Let
$\sigma\in\Gamma$.  By Lemmas~\ref{isomnormalizer} and~\ref{normalizer}, $\varphi_{\sigma}\colon X\to
\el{\sigma}{X}$ is given by $(M,e)$ where $\overline{M}\in G_{S_4}$.  A computation shows that $\sigma^*(t)$ is
one of the following:
\begin{enumerate}
  \item [i.] $t$
  \item [ii.] $-t$
  \item [iii.] $\frac{2t+12}{t-2} $
  \item [iv.] $\frac{2t-12}{-t-2}$
  \item [v.] $\frac{2t-12}{t+2}$
  \item [vi.] $\frac{2t+12}{-t+2}$.
\end{enumerate}

Since $\overline{\varphi}_{\tau} \colon  X/\Aut(X) \to \el{\tau}{(X/\Aut(X))}$ is defined over $K$ for all
$\tau\in\Gamma$, we have
$\overline{\varphi_{\tau}}\,\overline{\varphi_{\sigma}}=\overline{\varphi_{\sigma\tau}}$ for all
$\tau\in\Gamma$.  The fractional linear transformations i through vi form a group under composition isomorphic
to $S_3$. The map $\tau \mapsto \tau^*(t)$ defines a homomorphism from $\Gamma$ to this group. The kernel of
this homomorphism is $\Lambda:=\{\tau\in\Gamma\mid \tau^*(t)=t\}$.  So $|\Gamma/\Lambda|= 1, 2, 3, $ or $6$.
\begin{enumerate}
  \item [Case 1:] $|\Gamma/\Lambda|=1$.  In this case the fixed field of
$\Gamma^*$ is $K(t)$ and
  $B_K=\mathbb{P}^1_K$.
  \item [Case 2:] $|\Gamma/\Lambda|=2$.  Let $\sigma$ be a representative of the
  nontrivial coset.  There are three cases.
\begin{enumerate}
  \item  $\sigma^*(t)=-t$.  Then $t=0$ corresponds to a point $P\in B_K(K)$.
  \item $\sigma^*(t)=\frac{2t+12}{t-2}$.  Then $t=6$ corresponds to a point $P\in B_K(K)$.
  \item $\sigma^*(t)=\frac{2t-12}{-t-2}$.  Then $t=-6$ corresponds to a point ${P\in B_K(K)}$.
\end{enumerate}
 \item [Case 3:] $|\Gamma/\Lambda|=3$.  Since the fixed field of
$\Lambda^*$ is
  $\overline{K}^{\Lambda}(t)$, $B_K$ has a $\Kbar^\Lambda$-rational point.
By Lemma~\ref{odd}, since $[\overline{K}^{\Lambda}:K]$ is odd, $B_K$ has a $K$-rational point.
  \item [Case 4:] $|\Gamma/\Lambda|=6$.  Let $\Pi$ be a
  subgroup of $\Gamma$ containing $\Lambda$ such that
$\Pi/\Lambda$ is a subgroup of $\Gamma/\Lambda$ of order
  $2$. By
Case 2, $B_K$ has a $\ol K^{\Pi}$ rational point. Since $[\overline{K}^{\Pi}:K]=3$ is odd, by Lemma~\ref{odd},
 $B_K$ has a $K$-rational point.

\end{enumerate}

\item[(d)]  $G\cong A_4$.  The element $t':=x^2+x^{-2}$ is fixed by the normal subgroup $G_{V_4}$.  From (c), we
see that the element
\[t:=\frac{1}{4}t'\left(\frac{2t'-12}{t'+2}\right)\left(\frac{2t'+12}{-t'+2}\right)
=\frac{x^{12} - 33x^8 - 33x^4 + 1}{-x^{10} + 2x^6 - x^2}\] is fixed by $G_{A_4}$ and is a rational function of
degree $12$ in $x$.  So the function field of $X/\Aut(X)$ equals $\Kbar(t)$.  We use $t$ as coordinate on
$X/\Aut(X)$.  The map $\rho\colon X\rightarrow X/\Aut(X)$ is given by
\[(x,y)\mapsto (x^{12} - 33x^8 - 33x^4 + 1)/(-x^{10} + 2x^6 - x^2).\]Let $\sigma\in\Gamma$. By Lemmas~\ref{isomnormalizer} and~\ref{normalizer},
$\varphi_{\sigma}\colon X\to \el{\sigma}{X}$ is given by $(M,e)$ where $\overline{M}\in G_{S_4}$.    A
computation shows that $\sigma^*(t)=\pm t$.  Then $t=0$ corresponds to a point $P\in B_K(K)$.

\item[(e)]  $G\cong S_4$.  By Lemma~\ref{normalizer}, $N(G)=G$.  So by Corollary~\ref{def}, $X$ can be defined
over $K$.
 \item[(f)] $G\cong A_5$. By Lemma~\ref{normalizer}, $N(G)=G$.  So by Corollary~\ref{def}, $X$ can be
defined over $K$.

\item[(g)]  $G=G_{\beta, A}$. Let $d$ be the order of $\beta$ and let $t=g(x):=\prod_{\alpha\in A}(x-\alpha)^d$.
Then $t$ is a rational function of degree $|G|$ fixed by $G_{\beta, A}$ acting by fractional linear
transformations.  So the function field of $X/\Aut(X)$ equals $\ol K(t)$.  We use $t$ as a coordinate function
of $X/\Aut(X)$.  Let $\sigma\in\Gamma$. By Lemma~\ref{isomnormalizer}, $\varphi_\sigma\colon X\to
\el{\sigma}{X}$ is given by $(M,e)$ where $M$ is an upper diagonal matrix.  So $\sigma^*(t)=g^{\sigma}(ax+b)$
for some $a\ne 0$ and $b$.  Let $P$ be the point of $X/\Aut(X)$ corresponding to $x=\infty$. Then since
$g^{\sigma}(a\infty+b)=g(\infty)$, $P$ corresponds to a point in $B_K(K)$.

\item[(h)] $G=\PSL_2(\mathbb{F}_{p^r})$.  Let $q=p^r$.   It can be deduced from Theorem 6.21 on page 409
of~\cite{Suzuki:1} that $\PSL_2(\mathbb{F}_q)$ is generated by the image in $\PGL_2(\ol K)$ of the following
matrices
\[\left\{\begin{mat2}
           0&-1\\
           1&0\\
           \end{mat2},\begin{mat2}
           1&a\\
           0&1\\
           \end{mat2}\colon a\in\F_{p^r}\right\}.\]
Let \[g(x)=\frac{((x^q -x)^{q-1}+1)^{\frac{q+1}{2}}}{(x^q-x)^{\frac{q^2-q}{2}} }.\] One can verify that
$g(\frac{-1}{x})=g(x)$ and $g(x+a)=g(x)$ for all $a\in\F_{p^r}$.  Since $g$ is a rational function of $x$ of
degree $\frac{q^3-q}{2}=|\PSL_2(\F_q)|$, the function field of $X/\Aut(X)$ is $\ol K(t)$ where $t=g(x)$.  We use
$t$ as a coordinate function on $X/\Aut(X)$.  The map $\rho\colon X\rightarrow X/\Aut(X)$ is given by
\[(x,y)\mapsto \frac{((x^q
-x)^{q-1}+1)^{\frac{q+1}{2}}}{(x^q-x)^{\frac{q^2-q}{2}} }.\] Let $\sigma\in\Gamma$. By
Lemmas~\ref{isomnormalizer} and~\ref{normalizer}, $\varphi_{\sigma}\colon X\to \el{\sigma}{X}$ is given by
$(M,e)$ where $\overline{M}\in \PGL_2(\mathbb{F}_q)$.   A computation shows that $\sigma^*(t)=\pm t$.  Then
$t=0$ corresponds to a point $P\in B_K(K)$.

 \item[(i)] $G= \PGL_2(\mathbb{F}_{p^r})$.  By Lemma~\ref{normalizer}, $N(G)=G$.  So by
Corollary~\ref{def}, $X$ can be defined over $K$.

\end{enumerate}
\end{proof}

 Specific examples of hyperelliptic curves
not definable over their field of moduli are given on page 177 of~\cite{sh:1}; these examples have $|G|=1$.
Adjusting these examples, we now construct others with $|G|>5$.

Let $n>5$, let $m$ be odd, and consider the polynomial $f(x)\in\C[x]$ given by
\[f(x)=a_0x^{nm} +\sum_{r=1}^m(a_rx^{n(m+r)} + (-1)^ra_r^cx^{n(m-r)}),\]
with $a_m=1$, $a_0\in\mathbb{R}^*$, and where $z^c$ is the complex conjugate of $z$ for any $z\in\C$. Assume
that for $r=1,\ldots, m-1$ we have $a_r\ne(-1)^r\beta^{-nr}a_r^c$ for any $2mn^{th}$ root of unity $\beta$ and
that $f(x)$ is square free.
\begin{lem}\label{cyclic}Following the above notation, let $X$ by the hyperelliptic curve over $\C$ given by $y^2=f(x)$.  Let $\iota$ be the
hyperelliptic involution of $X$ and let $\nu$ be the automorphism of $X$ defined by ${\nu(x,y)=(\zeta x, y)}$,
where $\zeta$ is a primitive $n^{th}$ root of unity.  Then $\Aut(X)=\langle\iota\rangle\oplus\langle\nu\rangle$.
\end{lem}
\begin{proof} Let $G=\Aut(X)/\langle\iota\rangle$.  The image of $\nu$ in $G$ under the quotient map $\Aut(X)\to G$
has order $n$. Since $n>5$, by Lemma~\ref{subgroups}, $G$ is either cyclic or dihedral.  In either case the
image of $\nu$ in $G$ generates a cyclic normal subgroup of $G$.

Suppose that $G$ is cyclic of order $n'>n$.  Since the only elements in $\PGL_2(\C)$ that commute with the image
of diagonal matrices are the images of diagonal matrices, by Lemma~\ref{isom} there exists an element
$u\in\Aut(X)$ defined by
\[u(x,y)=(\zeta'x,ey)\]
where $e\in\C^*$ and $\zeta'$ is a primitive $(n')^{th}$ root of unity.  It follows that $f(\zeta' x)$ is a
scalar multiple of $f(x)$.  This is a contradiction by our choice of coefficients for $f$.

Suppose that $G$ is dihedral.  By Lemma~\ref{normalizer}~(a) and Lemma~\ref{isom}, there exists an element
$v\in\Aut(X)$ defined by
\[v(x,y)=(\alpha/x, e'y/x^{mn})\]
where $e'$, $\alpha\in\C^*$. It follows that $x^{2mn}f(\alpha/x)$ is a scalar multiple of $f(x)$. Since
\[x^{2mn}f(\alpha/x)=\alpha^{nm}(a_0x^{nm}+\sum_{r=1}^m((-1)^r\alpha^{-nr}a_r^cx^{n(m+r)}+\alpha^ra_rx^{n(m-r)})\]
and $a_0\ne 0$, we must have $\frac{x^{2mn}}{\alpha^{nm}}f(\alpha/x)=f(x)$.  Since $a_m=1$, we must have
$\alpha^{mn}=-1$ and $a_r=(-1)^r\alpha^{-nr}a_r^c$ for $r=1,\ldots, m-1$. This is a contradiction. Therefore $G$
is cyclic of order $n$.

The function field of $X$ is $\C(x,y)$ and the function field of $X/\Aut(X)$ is $\C(x^n)$.  Since the places in
$\C(x^n)$ corresponding to $x^n=0$ and ${x^n=\infty}$ do not ramify completely in $\C(x,y)$, by Theorem~5.1
of~\cite{Brandt:1} we have ${\Aut(X)=\langle\iota\rangle\oplus\langle\nu\rangle}$.
\end{proof}
\begin{prop}\label{notdef}Following the above notation, let $X$ by the hyperelliptic curve of genus $g=mn-1$
over $\C$ given by $y^2=f(x)$. The field of moduli of $X$ relative to the extension $\C/\R$ is $\R$ and is not a
field of definition for $X$.
\end{prop}
\begin{proof}  By Lemma~\ref{cyclic}, $\Aut(X)=\langle\iota\rangle\oplus\langle\nu\rangle$ where $\iota$ is
the hyperelliptic involution of $X$, and $\nu(x,y)=(\zeta x, y)$ where $\zeta$ is a primitive $n^{th}$ root of
unity.  The map $\mu$ defined by
\[\mu(x,y)=((\omega x)^{-1},ix^{-nm}y),\]
where $\omega^n=-1$, is an isomorphism between the curve $X$ and the complex conjugate curve $\el{c}{X}$. Any
isomorphism $X\to \el{c}{X}$ is given by $\mu\nu^k$, or $\mu \iota\nu^k$ for some $0\le k\le n-1$.  We have $\mu
\iota=\iota\mu$,
 \[ \mu\nu(x,y)=((\omega\zeta x)^{-1},i(\zeta x)^{-nm}y)=\nu^c\mu(x,y),\] and
 \[\mu^c\mu(x,y)=( (\omega^{-1}(\omega x)^{-1})^{-1},-i(\omega x)^{nm}(ix^{-nm}y))=(\omega^2 x, -y)=\nu^l \iota(x,y)\]
for some $l$. Then
\[(\mu\nu^k)^c\mu\nu^k=\mu^c\nu^{-k}\mu\nu^k=\mu^c\mu\nu^{2k}=\iota\nu^{2k+l}\ne Id\]
and
\[(\mu \iota\nu^k)^c\mu \iota\nu^k=\mu^c \iota\nu^{-k}\mu \iota\nu^k=\mu^c\mu\nu^{2k}=\iota\nu^{2k+l}\ne Id.\]
Therefore Weil's cocycle condition from Theorem~1 of~\cite{Weil:1} does not hold.  So $X$ cannot be defined over
$\mathbb{R}$.
\end{proof}

\section*{Acknowledgements}
The author learned of the conjecture of~\cite{Shaska:1} from Bjorn Poonen, and thanks him for comments on an
early draft of this paper. The author was partially supported by NSF grant DMS-0301280 of Bjorn Poonen.

\bibliographystyle{amsplain}

\bibliography{fom}
\small{Department of Mathematics, University of California, Berkeley, CA 94720
\newline E-mail address: bhuggins@math.berkeley.edu}

\end{document}